\documentclass[12pt]{elsarticle}

\usepackage[utf8]{inputenc} 

\usepackage[margin=1.in]{geometry} 
\usepackage{graphicx} 
\usepackage[dvipsnames]{xcolor}

\usepackage{array} 
\usepackage{verbatim} 

\usepackage{amsthm,amsmath,amssymb}

\usepackage[english]{babel}
\usepackage{graphicx}

\usepackage[small,bf,hang]{caption} 

\usepackage{subfig}

\usepackage{hyperref}

\journal{???}

 \graphicspath{{figures/}}

\begin{document}
\begin{frontmatter}

\title{House of Graphs 2.0: a database of interesting graphs and more}

\author[gent]{Kris Coolsaet}
\ead{kris.coolsaet@ugent.be}

\author[gent,kulak]{Sven D'hondt}
\ead{dhondtsven36@gmail.com}

\author[gent,kulak]{Jan Goedgebeur\fnref{kul}}
\ead{jan.goedgebeur@kuleuven.be}

\address[gent]{Department of Applied Mathematics, Computer Science \& Statistics,\\
  Ghent University, 9000 Ghent, Belgium\medskip}

\address[kulak]{Department of Computer Science,\\
  KU Leuven Kulak, 8500 Kortrijk, Belgium}

\fntext[kul]{Supported by Internal Funds of KU Leuven.}

\begin{abstract}
In 2012 we announced ``the House of Graphs'' (\url{https://houseofgraphs.org}) [Discrete Appl.\ Math.\ 161 (2013), 311-314], which was a new database of graphs. The House of Graphs hosts complete lists of graphs of various graph classes, but its main feature is a searchable database of so called ``interesting'' graphs, which includes graphs that already occurred as extremal graphs or as counterexamples to conjectures. An important aspect of this database is that it can be extended by users of the website.

Over the years, several new features and graph invariants were added to the House of Graphs and users uploaded many interesting graphs to the website. But as the development of the original House of Graphs website started in 2010, the underlying frameworks and technologies of the website became outdated. This is why we completely rebuilt the House of Graphs using modern frameworks to build a maintainable and expandable web application that is future-proof. On top of this, several new functionalities were added to improve the application and the user experience.

This article describes the changes and new features of the new House of Graphs website.
\end{abstract}


\end{frontmatter}


\section{Introduction}
\label{sect:introduction}

The website \textit{the House of Graphs} (HoG) (\url{https://houseofgraphs.org}) -- which we originally announced in 2012 in~\cite{HoG} -- basically consists of two parts. In the ``graph meta-directory'' we host complete lists of all pairwise non-isomorphic graphs up to a given order of various graph classes (e.g.\ cubic graphs, fullerenes, planar graphs, snarks, trees, etc.) together with links to the papers describing the algorithms which were used to generate these graphs and their source code. But its main part is the ``searchable database of interesting graphs''.

Most graph theorists will agree that among the vast number of graphs that exist, there are only a few that can be considered really ``interesting''. So as we already paraphrased Orwell's famous words in~\cite{HoG}, we could say that: \textit{``All graphs are interesting, but some graphs are more interesting than others''}.
We felt that it would be useful to have a database which allows users to search for graphs which are particularly relevant for the problem or conjecture they are working on (e.g.\ to gain more intuition about the problem or to find a counterexample).

The main reason why we decided to start the development of the original House of Graphs website in 2010, is that we wanted to fill this gap and provide a useful service to the graph theory community. Moreover, we also wanted to serve as a central repository for complete lists of graphs of various graph classes as previously such lists were scattered over several websites hosted by individual researchers (see~\cite{HoG} for concrete examples).

Note that nor in the previous article~\cite{HoG}, nor here we will try to give an exact definition of which graphs are ``interesting'' as this depends on the problem one wants to study. Graphs which we definitely consider as ``interesting'' include named graphs in the literature (e.g.\ the Coxeter graph, the Heawood graph, the Petersen graph, etc.), counterexamples to known conjectures or extremal graphs pointed out by a conjecture-making system. 

For every ``interesting'' graph in the searchable database, several graph invariants (e.g.\ chromatic number, diameter, girth,... -- see Table~\ref{table:invariants} for a complete list) are precomputed and stored in the database. Users of the web application can enter queries (e.g.\ search for graphs with chromatic number 4 and girth 5) and browse all graphs in the database matching the query (which can still be refined later). One can then inspect these graphs, study the values of the invariants or download the graphs in different formats for further personal use.
A very important functionality of the website is that users can also upload new graphs which they find interesting. This way, the database of interesting graphs keeps growing to become as complete as possible.

Thus the main functionality of the House of Graphs is to offer a searchable database of ``interesting'' graphs which yields a relatively small list of graphs matching the users' queries that still gives them a good chance to find counterexamples or obtain results that allow to judge
the general situation. 

When the House of Graphs website was launched in 2012, the searchable database was initialised with about 1\,500 ``interesting'' graphs. Most of these graphs were extremal graphs pointed out by the conjecture-making system \textit{GraPHedron}~\cite{melot2008facet} (see~\cite{HoG} for details) and it also included graphs which occurred as a counterexample for certain conjectures (see e.g.~\cite{BGHM13}). In the years that followed, many new ``interesting'' graphs were added to the database and currently it contains nearly 22\,000 graphs. (Note that we intentionally keep the searchable database relatively small in size. This is because we want the number of graphs matching a query to be manageable for inspection by the users and only to contain the graphs that are likely to be relevant for the users.)

Because the development of the House of Graphs started in 2010, the underlying frameworks and technologies became outdated over the years. However, the content of the database and the need to use these graphs for research is still relevant and will be relevant for the foreseeable future. This is why in 2021--2022, the House of Graphs was rebuilt completely, using modern frameworks to build a maintainable and expandable web application that is future-proof. On top of this, new functionalities were added to improve the application and the user experience. 

One of these improvements focuses on calculating graph invariants using stand-alone programs, not tied to a particular programming language. In the old House of Graphs these invariants were calculated using \textit{Grinvin} (a Java graph library which is described in~\cite{peeters2009grinvin}), which significantly limited the possibilities for expansion and adaptation. The new flexible invariant computation framework allowed us to add several new graph invariants (such as: automorphism group size, number of vertex orbits, treewidth, traceability, hypohamiltonicity,...) and to replace several existing invariant computers by faster implementations (see Section~\ref{subsect:invariants} for more details).

Next to this, a new tool for drawing and editing graphs was developed which also allows to import graphs from a file to generate an automated drawing which can then be altered manually. It is now also possible to edit the drawings of a given graph and add multiple drawings per graph. Moreover, the graph drawings can be exported to various formats (svg, png, pdf, tikz,...). Furthermore, the new website renders much more nicely on mobile devices.

This article is organised as follows. In Section~\ref{sect:technologies} we give an overview of the technologies and frameworks used in the new House of Graphs website.
In Section~\ref{sect:functionality} we describe the functionality of the website into more detail, in particular: the graph meta-directory (Section~\ref{subsect:meta-dir}), the searchable database (Section~\ref{subsect:searchable_db}), and the graph drawing utility (Section~\ref{subsect:drawing}). Finally, we conclude this article in Section~\ref{sect:conclusion} with an outlook on possible directions for future work.


\section{Technologies}
\label{sect:technologies}

When we look at the old setup of the House of Graphs, two distinguishable parts can be extracted from the application. The first part of the application is the database, which uses \textit{PostgreSQL}~\cite{PostgreSQL} as its relational database management system. The rest of the application was managed by the \textit{Apache Struts 2.3}~\cite{ApacheStruts} framework. This is an open-source framework for developing Java EE web applications. It encourages developers to adopt the model-view-controller architecture and uses Java Server Pages to display the web pages. But over the years the 2.3 version of this framework has reached end-of-life status and it would have required a lot of work to migrate it to the most recent version of the framework. Instead, it seemed better to us to completely rebuild the new website (if possible using a more suitable and more future-proof framework) and seize the opportunity to add several new features and enhancements which improve the user experience.

\begin{figure}[h!tb]
	\centering
	\includegraphics[width=0.55\textwidth]{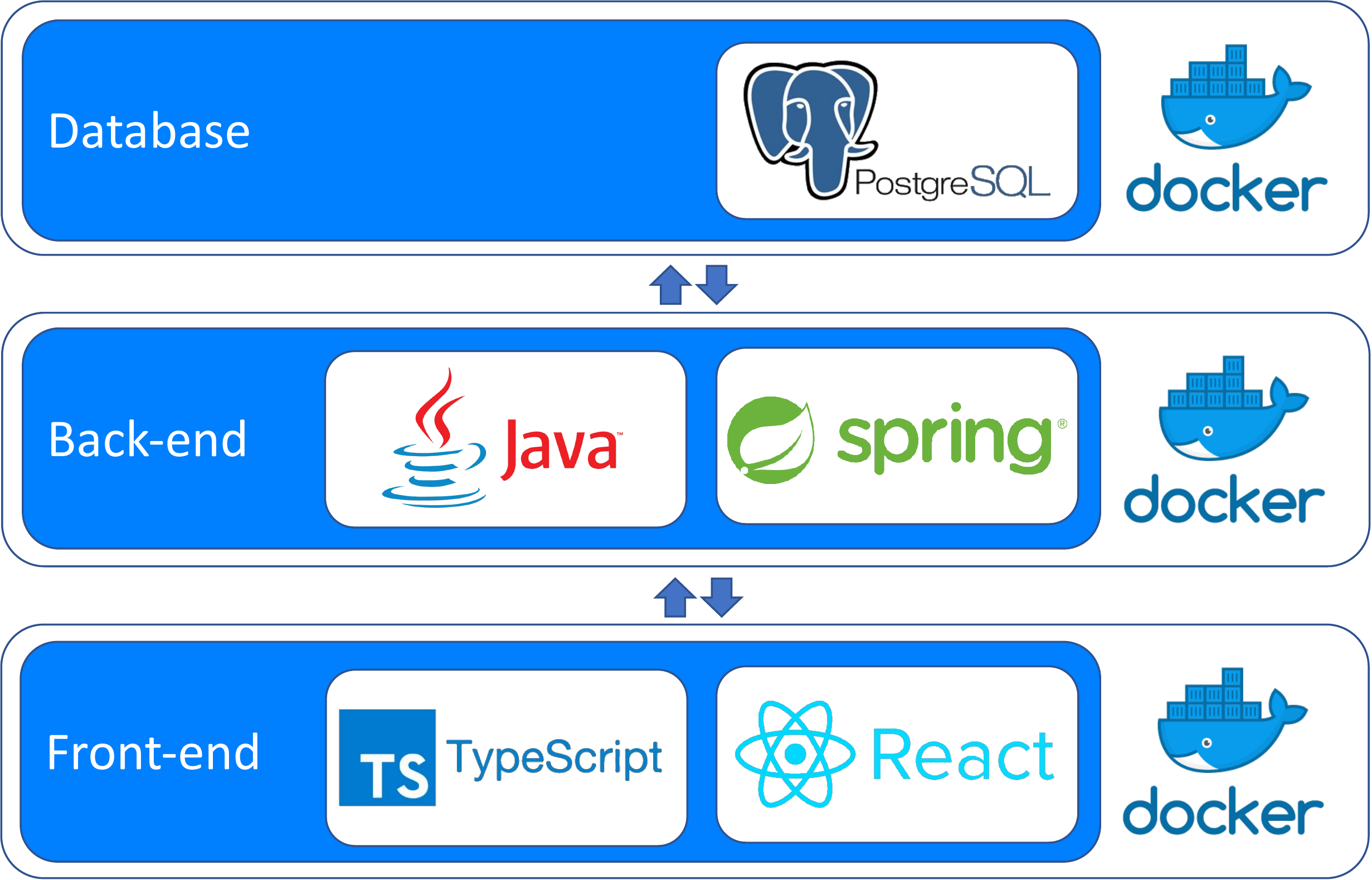}
	\caption{An overview of the frameworks used in the new House of Graphs website.}
	\label{fig:technologies}
\end{figure}

To construct a new web application for the House of Graphs, the most important aspect was to focus on modern, popular frameworks that are future-proof. Modern web application stacks often separate the application in three main parts, being the database, the back-end, and front-end. As \textit{PostgreSQL} is still very popular today, this framework was also used for the new database. This also guarantees a smooth transition of the data. To build the front- and back-end, new technologies and frameworks were used. Although \textit{Struts} is still being maintained, other frameworks have gained much more success over the years. The \textit{Java Spring Framework}~\cite{Spring} was selected as framework to construct the back-end of the new application, while \textit{React}~\cite{React} with TypeScript was used for the front-end. Both frameworks are widely used and are guaranteed to continue development and maintenance in the foreseeable future. Next to these frameworks, \textit{Maven}~\cite{Maven} and \textit{Docker}~\cite{Docker} were used for automatically building, executing and deploying the application. Docker allows to run the applications of the different layers (database, back-end and front-end) in different containers so they can be easily migrated to a new server or distributed over different servers in the future, if needed. Figure~\ref{fig:technologies} gives a schematic overview of the technologies used and the interaction between the layers.

\subsection{Migrating the database}

As stated in Section~\ref{sect:technologies}, the old House of Graphs used a \textit{PostgreSQL} database. The interaction with this database was managed by the \textit{Hibernate ORM} framework~\cite{Hibernate}. This is an open-source Java framework that provides an object-relational mapping between the object-oriented domain model and the relational database.

As \textit{PostgreSQL} is still supported and widely used today, there was no reason to change the database framework. However, to support the functionalities of the new web application, the old database schema needed to be altered. One of the new functionalities of the new House of Graphs is to store multiple graph embeddings (visualisations) with a single graph. This had to be taken into account while designing the new database schema, as a separate table for embeddings is now needed. To accommodate for safe authentication, an extra table was required for storing activation tokens. Due to the use of \textit{Hibernate}, some unnecessary tables were present in the database that only establish a join relationship between two other tables. These tables can therefore safely be left out. More details about the old and new database schema can be found in~\cite{thesissven}.

\textit{Hibernate ORM} has proven to come with some disadvantages in the old application, as \textit{Hibernate} functions as an abstraction layer. As a consequence, the developer is shielded from what happens at the level of the database, which makes it harder to intervene. Due to this abstraction, some knowledge of the framework is required to really understand what is happening behind the scenes. \textit{Hibernate} also introduces some overhead when used with small database schemas. Therefore, we prefer a framework that lies closer to the database itself. This is for example useful to manually construct specific search queries for graphs (e.g.\ using formulas involving invariants). Instead of \textit{Hibernate}, the new House of Graphs application uses \textit{DAOHelper}, a Java library which was developed by the first author. This library provides a fluent interface for creating SQL calls that are closely related to \textit{PostgreSQL} queries without introducing overhead.

With the new database schema in place, we established a migration process to migrate the data from the original database to the new database, so all graphs which were uploaded to the old website as well as the comments which were associated with the graphs are preserved. We also migrated the accounts of the registered users, but for security reasons the passwords were reset as we now use a stronger algorithm to hash the passwords and increased the minimum length of the passwords.

In particular, previously we used MD5 as a hashing algorithm. However, it has been shown that MD5 is not collision resistant~\cite{md5_break}. As a consequence, methods exist to detect hash collisions that perform better than brute-force techniques. On top of that, no salt was used to store the passwords. This salt is unique for each user and is added to the password before hashing it. This makes sure that hashing the same passwords gives different results, which significantly hampers attackers. In the new web application, we therefore opted for \textit{bcrypt}~\cite{bcrypt}, which is also recommended by The Open Web Application Security Project (OSWAP)~\cite{owasp}. \textit{Bcrypt} is collision resistant and incorporates a salt in its hashing algorithm. 

We also made sure that the HoG graph id of the graphs which were present in the original database remained the same in the new database as researchers sometimes mention these id's in papers to refer to specific graphs on the House of Graphs. 

The new website also came with a new url: \url{https://houseofgraphs.org}. But the old url \url{https://hog.grinvin.org} also points to the new website and will remain active in the foreseeable future for backwards compatibility.

The url structure was also changed to make it more structured and logical. For example, the old url for the meta-directory page with all cubic graphs up to a given order was \url{https://hog.grinvin.org/Cubic}, but it now is \url{https://houseofgraphs.org/meta-directory/cubic}. Similarly, previously the direct url with the graph details of the graph with HoG graph id 26 was \url{https://hog.grinvin.org/ViewGraphInfo.action?id=26} while it now is \url{https://houseofgraphs.org/graphs/26}. However, we added redirection-mappings so these old urls will keep working as such direct links were mentioned in several previously published papers.


\section{Functionality of the website}
\label{sect:functionality}

\subsection{The graph meta-directory}
\label{subsect:meta-dir}

In the graph meta-directory we offer static lists of all pairwise non-isomorphic graphs up to a given order for various graph classes. Table~\ref{table:graph_classes} gives an overview of the graph classes that are currently listed in the meta-directory. In many cases we also list subclasses of a given graph class in the metadirectory. For example: the meta-directory page on snarks also contains separate more fine-grained sublists of flow-critical snarks, hypohamiltonian snarks, and permutation snarks.

\begin{table}[htb!]
\centering
\footnotesize
	\begin{tabular}{|l | l | l |}	
		\hline
\textit{Almost hypohamiltonian graphs} & \textit{Largest degree-diameter graphs*}  & \textit{Ramsey numbers} \\
\textit{Alternating plane graphs}  & Maximal triangle-free graphs & Regular graphs*\\
\textit{Cographs} & \textit{Minimal Cayley graphs} & Simple graphs*\\
\textit{Critical $H$-free graphs} & \textit{Minimal Ramsey graphs} & Snarks\\
Cubic graphs  & \textit{Nut graphs} & Strongly regular graphs*\\
\textit{Directed graphs*} & \textit{Perihamiltonian graphs}   & Trees\\
Fullerenes & Planar graphs  & \textit{Triangle-free $k$-chromatic graphs}\\
\textit{Hypohamiltonian graphs} & \textit{Platypus graphs}  &  \textit{Uniquely hamiltonian graphs}\\
\textit{Interval graphs*}  & \textit{Quartic graphs}  & Vertex-transitive graphs*\\
		\hline
	\end{tabular}
\caption{The graph classes which are currently listed in the graph meta-directory. The lists which are new since the announcement of the original House of Graphs website in~\cite{HoG} are in italic. The lists which are hosted elsewhere are marked with an asterisk.}
\label{table:graph_classes}
\end{table}

A lot of new lists were added compared to the old House of Graphs (i.e.\ in 2012 we had 10 and currently we offer 27 graph classes) and we plan to add more graph classes in the future.
Most of these lists are hosted on the House of Graphs, but some are direct links to pages hosted by other people. 

These complete lists of graphs can be downloaded in \textit{graph6} format up to the orders for which it is still feasible to store these graphs. Graph6 is a format introduced by McKay~\cite{graph6} that is commonly used to encode graphs. Many mathematical software packages such as \textit{SageMath}~\cite{sagemath} can handle graphs in this format.

Next to that, we also refer to the papers which describe the algorithms which were used to generate these graphs (e.g.\ \textit{buckygen}~\cite{BGM12}, \textit{geng}~\cite{nauty-website, MP14}, \textit{genreg}~\cite{M99}, and \textit{snarkhunter}~\cite{BGHM13}) and how the completeness of the lists was independently verified. In most cases these generators are open source and we link to their source code so researchers can directly generate these graphs themselves (also for the orders for which it was not possible to store all of them on the website).

\subsection{The searchable database of interesting graphs}
\label{subsect:searchable_db}

The main feature of the House of Graphs is its searchable database of ``interesting'' graphs. As explained in Section~\ref{sect:introduction}, we also offer a database of ``interesting'' graphs which includes named graphs in the literature, counterexamples to known conjectures, extremal graphs pointed out by a conjecture-making system and other ``interesting'' graphs which were uploaded by the users. Its main goal is to allow users to search for graphs that are relevant for the problem or conjecture they are working on in such a way that the number of graphs matching their query is relatively small so that they can be manually inspected but is still representative enough to allow to gain an intuition and to judge the general situation.

For each graph which is present in the database, several graph invariants were precomputed and stored in the database. In Section~\ref{subsect:invariants} we list these invariants and in Section~\ref{subsect:upload} we will explain how these computations are launched when uploading new graphs to the website.
As will be explained into more detail in Section~\ref{subsect:performing_search}, based on these invariants users can enter queries such as ``chromatic number $\geq 5$'' which in this specific example results in an overview of all ``interesting'' graphs with chromatic number at least 5 in the database. These graphs can then be inspected or exported by the user and these search queries can also be combined or refined by the user.

\subsubsection{Invariants}
\label{subsect:invariants}
Table~\ref{table:invariants} lists all boolean and numeric graph invariants which are currently supported by the website. When the original House of Graphs website was launched in 2012, it supported 24 graph invariants while the new website currently supports 44 invariants. On the old website we used the \textit{Grinvin} framework to compute the invariants. \textit{Grinvin} is a Java software package that was developed in~\cite{peeters2009grinvin} for the study of graphs and their invariants. It was specifically designed to teach students about graphs and mathematical reasoning. Although this approach of calculating invariants works correctly, there were some important limitations related to it which we wanted to handle in the new application. In particular, we wanted to disconnect the calculation process from the \textit{Grinvin} library as \textit{Grinvin} was not developed to be used as a calculation framework. As a consequence, the interaction with \textit{Grinvin} is not optimal and many invariants are not calculated in the most efficient way. A second disadvantage is that it was not possible to add new invariants to the House of Graphs that are not present in \textit{Grinvin}. To add a new invariant, the algorithm needs to be included in \textit{Grinvin} first (and must be implemented in Java using the predescribed \textit{Grinvin} interfaces and data format), after which it can be used by the House of Graphs. This process is very cumbersome and far from ideal.

\begin{table}[htbp!]
\scriptsize
\centering
\begin{tabular}{| l || l | l | l |}
\hline
Boolean invariants & \multicolumn{3}{c}{Numeric invariants} \vline\\
\hline
Acyclic & Algebraic Connectivity  & Girth & Number of Edges \\
Bipartite & Average Degree &  \textit{Group Size} & \textit{Number of Spanning Trees} \\
\textit{Claw-Free} & Chromatic Index  & Independence Number & \textit{Number of Triangles} \\
Connected & Chromatic Number  & Index & \textit{Number of Vertex Orbits} \\
\textit{Eulerian} & \textit{Circumference}  & Laplacian Largest Eigenvalue & Number of Vertices \\
\textit{Hamiltonian} & Clique Number  & \textit{Longest Induced Cycle} & \textit{Number of Zero Eigenvalues} \\
\textit{Hypohamiltonian} & \textit{Density}  & \textit{Longest Induced Path} & Radius \\
\textit{Hypotraceable} & Diameter  & Matching Number & Second Largest Eigenvalue \\
\textit{Planar} & Domination Number  & Maximum Degree & Smallest Eigenvalue \\
Regular & Edge Connectivity  & Minimum Degree & \textit{Treewidth} \\
\textit{Traceable} & \textit{Genus}  & \textit{Number of Components} & \textit{Vertex Connectivity} \\
\hline
\end{tabular}
\caption{The boolean and numeric graph invariants which are currently supported. The invariants which are new since the announcement of the original House of Graphs website in~\cite{HoG} are in italic.}
\label{table:invariants}
\end{table}


This is why in the new House of Graphs, we decoupled this library so it is now possible to calculate the invariants using stand-alone programs that can be called from the command line. Depending on the program itself, an input file can then be fed to the program, after which the output of the program can be read and stored in the database. The main advantage of this approach is that it becomes a lot easier to add new invariants. A second advantage is that it now is very simple to change the algorithm of a certain invariant. If one of the \textit{Grinvin} implementations has a bad performance, a new program can be created and can simply be plugged in to replace the existing implementation. An additional advantage is that this program can be written in any programming language. This means that a slow Java implementation can be replaced by a faster implementation that is written in C. A final advantage is that we can have multiple programs to compute a given invariant so that this could be used to independently verify these results.

We replaced several slow Java implementations of invariant computers from \textit{Grinvin} by much faster algorithms (e.g.\ for hamiltonicity, 
connectivity, independence number) and added several new invariants which were not present in \textit{Grinvin} (e.g.\ planarity, genus, treewidth,...). Some of these new invariant computers were implemented by us in the course of our research, while others were implementations of other researchers which were available as open source software. A selection of the invariant computers we used include \textit{nauty}~\cite{nauty-website, MP14}, the planarity algorithm from~\cite{boyer2004cutting}, the genus tester from~\cite{brinkmann2022practical}, the (hypo)hamiltonicity and (hypo)traceability tester from~\cite{K2hypo}, and the treewidth tester from~\cite{tamaki2019positive}.

\subsubsection{Performing searches and inspecting the results}
\label{subsect:performing_search}

The database can be searched in various ways:

\begin{itemize}
\item It is possible to search for graphs where a given numeric invariant has a certain value or is within a certain range (e.g.\ graphs where the independence number is between 5 and 7).
\item For boolean invariants (cf.\ Table~\ref{table:invariants}) one can specify if they must be \verb|true| or \verb|false| (e.g.\ only search for non-hamiltonian graphs).
\item One can search for graphs which were marked as ``interesting'' for a specific invariant by the user who uploaded the graph. (E.g.\ when one is working on a conjecture related to the chromatic number, it would make sense to search for graphs which are marked as ``interesting'' for this invariant).
\item One can also search for graphs which are associated with a given text or keyword. This searches for graphs where the search string occurs in the graph name or in the comments associated with this graph. (E.g.\ when one searches for the text ``petersen'' the resulting list of graphs will include the Petersen graph, but also the Petersen Cone and Generalised Petersen graph $GP(7,2)$).
\end{itemize}

\begin{figure}[h!tb]
	\centering
	\includegraphics[width=0.8\textwidth]{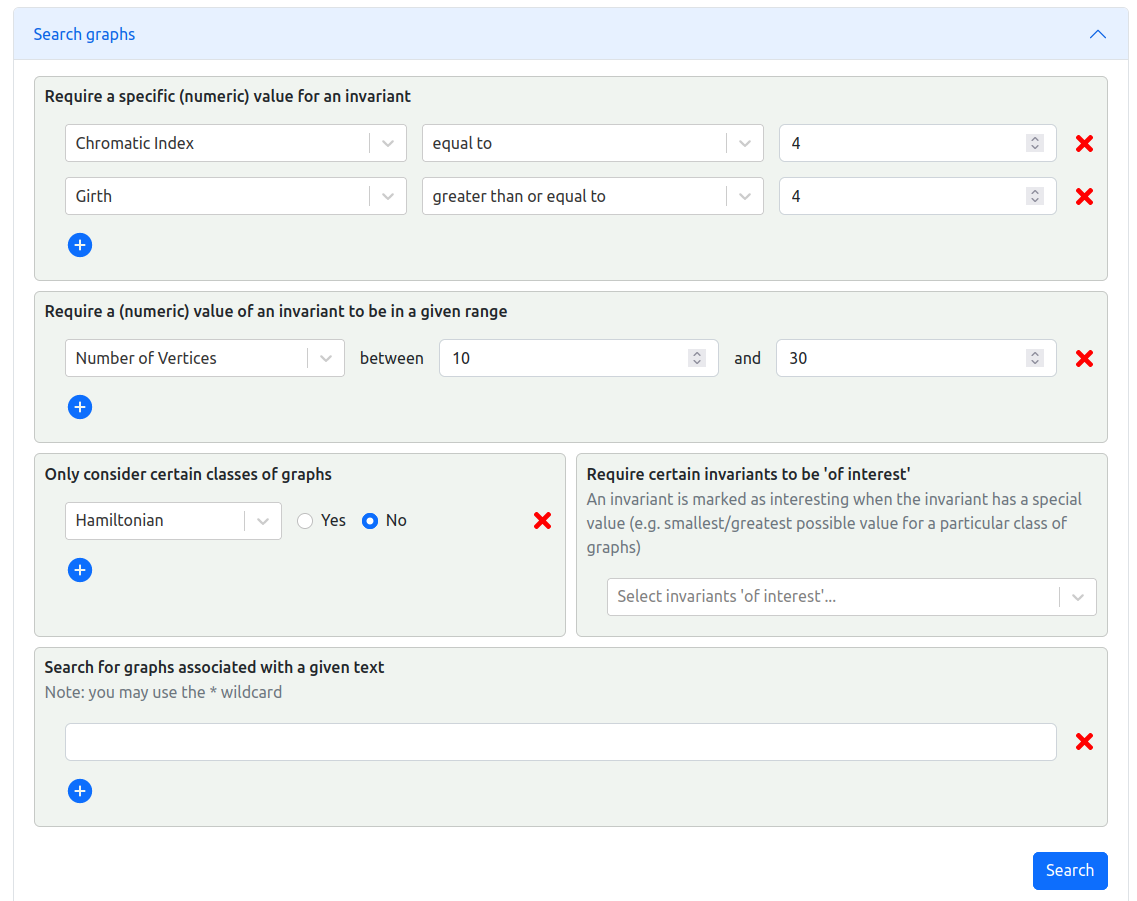}
	\caption{Screenshot of the page to enter the search queries.}
	\label{fig:search_graph}
\end{figure}

These search queries can also be combined (on the old House of Graphs website all criteria had to be entered one by one before the query could be refined in the next step). Figure~\ref{fig:search_graph} shows how these search queries can be entered on the website. In the next stage it is also possible to further refine the search (if wanted).

\begin{figure}[h!tb]
	\centering
\begin{minipage}[c]{.57\textwidth}
\centering   
   \subfloat[]{\label{fig:search_results_new}\includegraphics[width=1.0\textwidth]{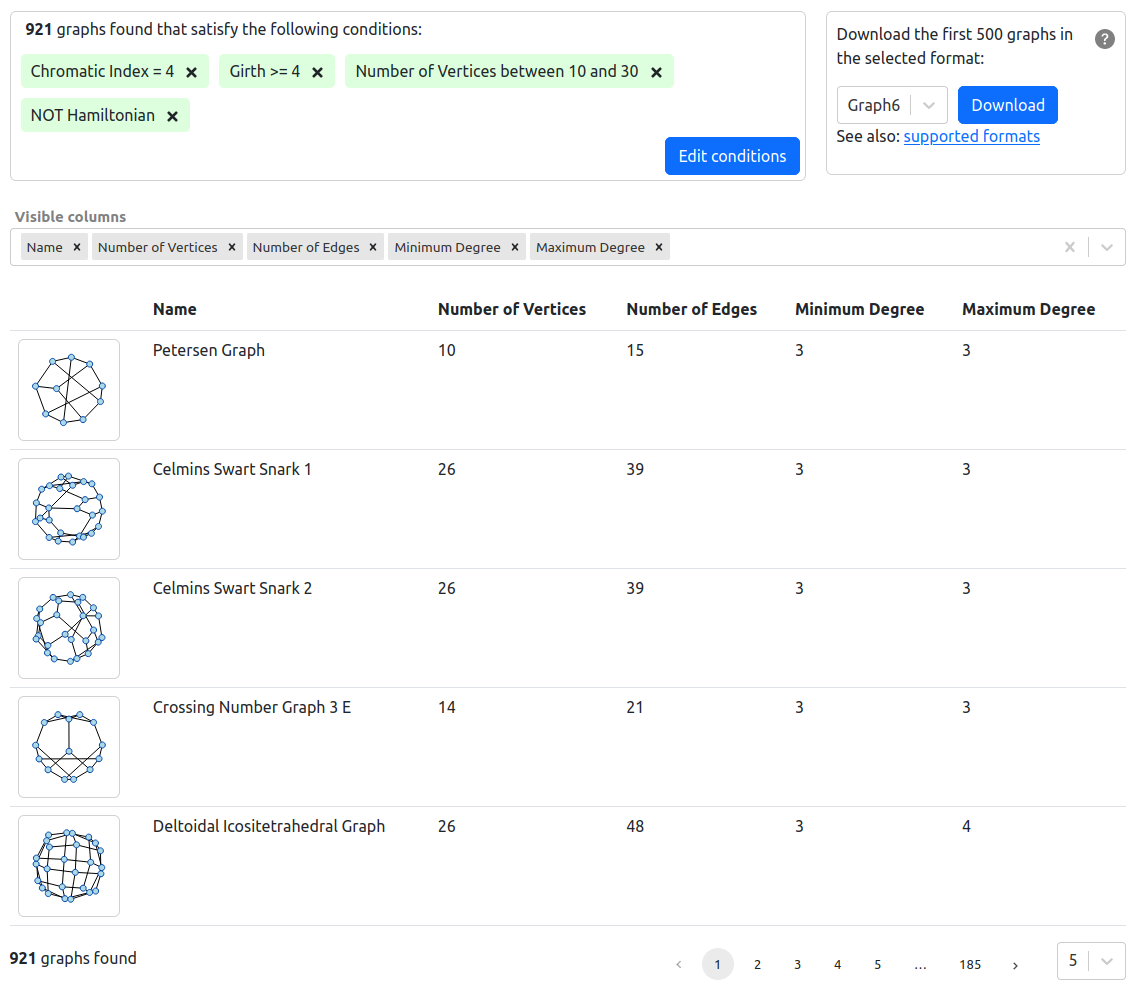}} 
\end{minipage}     	
\begin{minipage}[c]{.4\textwidth}
\centering   
   \subfloat[]{\label{fig:search_results_old}\includegraphics[width=0.8\textwidth]{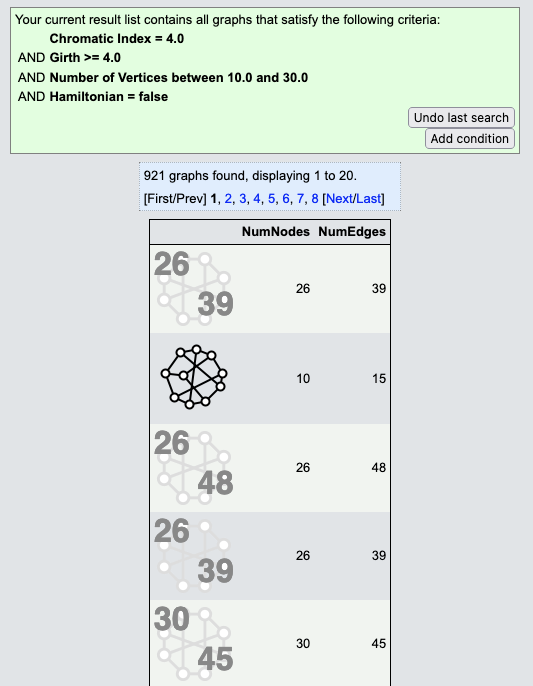}} 
\end{minipage}  	
	\caption{Figure~\ref{fig:search_results_new} is a screenshot of the page with the graphs matching the search query from Figure~\ref{fig:search_graph}. For comparison, Figure~\ref{fig:search_results_old} shows what this looked like on the old website.}
	\label{fig:search_results_both}	
\end{figure}

The website then gives an overview of all graphs in the database which match your search query. An example of such an overview is shown in Figure~\ref{fig:search_results_both}. The website shows a table where every row corresponds to a graph and the columns contain a drawing of the graph, the name of the graph (if applicable), and the values of the invariants which the user asked to see (of course this list of columns can be edited by the user). If the resulting list of graphs is too big, the user can still refine the search.

\begin{figure}[h!tb]
	\centering
	\includegraphics[width=0.8\textwidth]{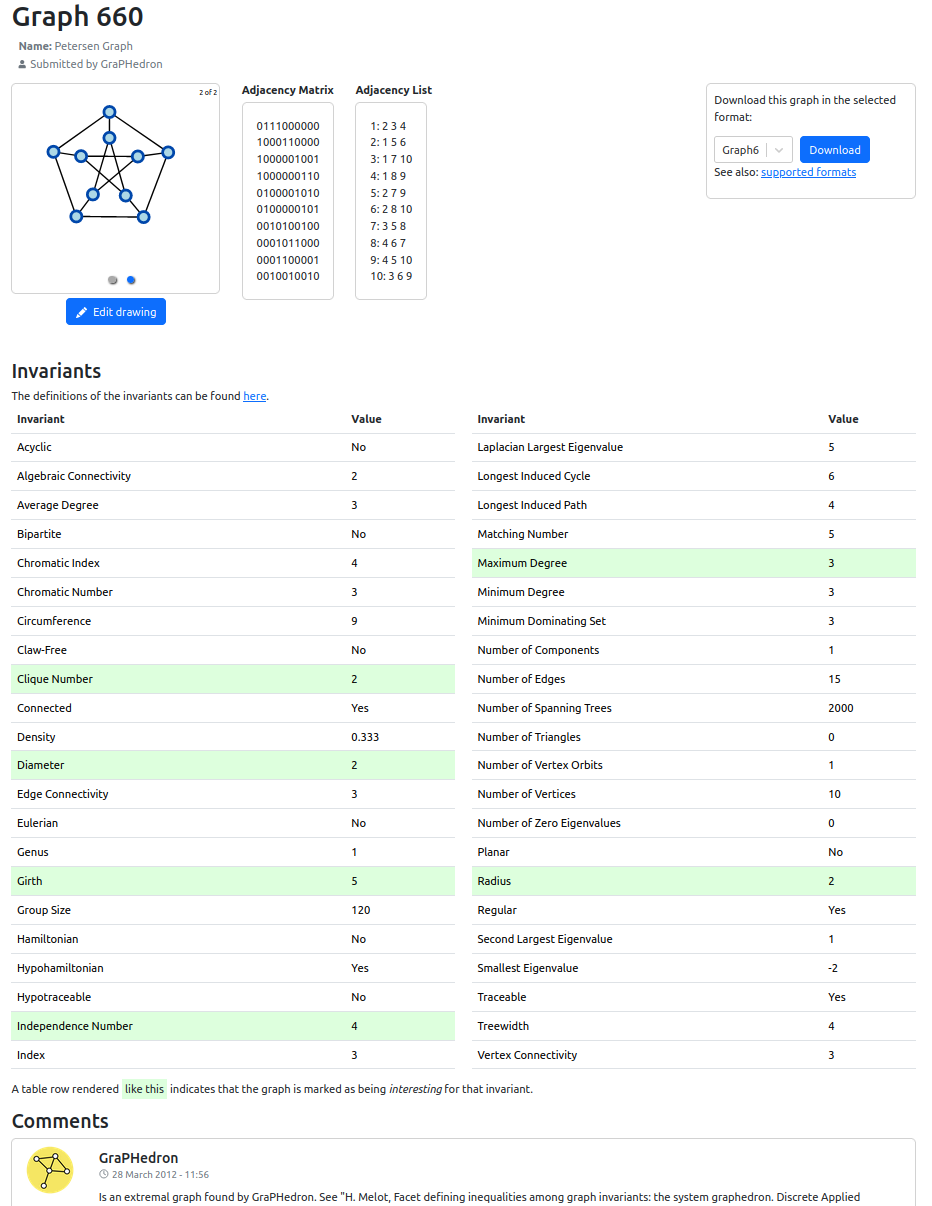}
	\caption{Screenshot of the page with the details of a given graph.}
	\label{fig:graph_details}
\end{figure}

\begin{figure}[h!tb]
	\centering
\begin{minipage}[c]{.3\textwidth}
\centering   
   \subfloat[]{\includegraphics[width=1.0\textwidth]{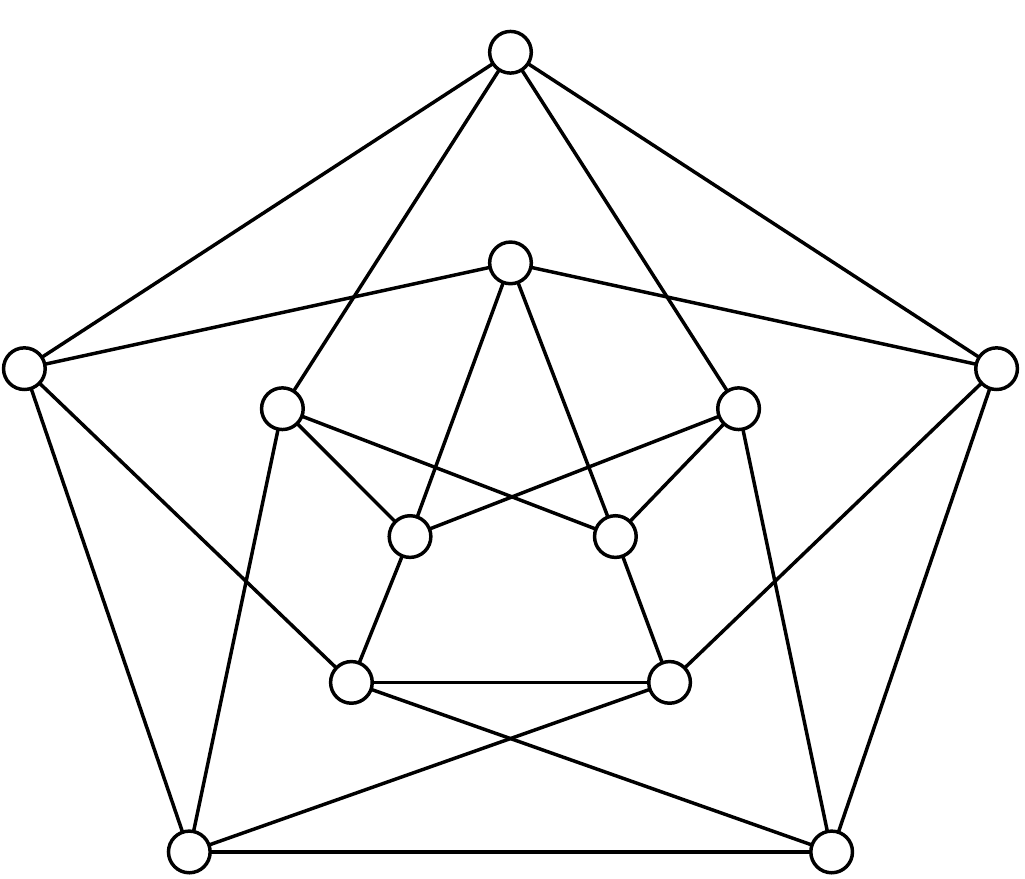}} 
\end{minipage}     	
\qquad \qquad
\begin{minipage}[c]{.3\textwidth}
\centering   
   \subfloat[]{\includegraphics[width=1.0\textwidth]{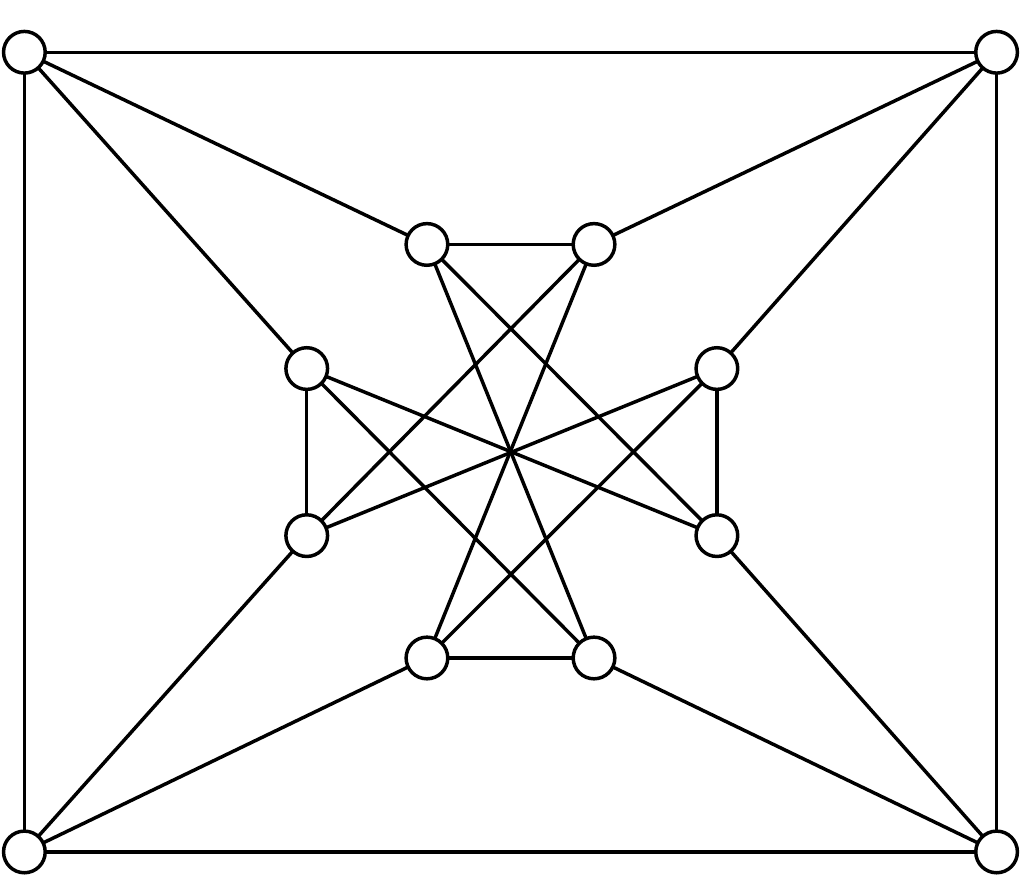}} 
\end{minipage}  	
	\caption{Two drawings of the Chv\'atal graph.}
	\label{fig:chvatal_graph}
\end{figure}

A user can also click on a specific graph in the result list and is then taken to a separate page which shows the details of this graph (see Figure~\ref{fig:graph_details}). In particular, this page shows the HoG graph id, the name of the graph (if applicable), who uploaded the graph, the adjacency matrix and adjacency list of the graph, all computed invariant values of this graph, and any comments associated with this graph (e.g.\ in what context it occurred). Moreover, this page also shows one or more drawings of this graph. This can be an automated drawing generated by the HoG spring embedder, but also drawings uploaded by the users using the graph editor (see Section~\ref{subsect:drawing}). Any user can open an existing embedding via the website and manipulate or export it to one of the supported formats (svg, png, pdf, tikz,...), but to avoid abuse only registered users can add new drawings to a graph (and similarly only registered users can add additional comments to a graph). This is an important new feature of the website as previously it was only possible to store one embedding (which could not be modified later), while there are usually many ``nice'' ways to draw a given graph (depending on the context one is working in). For example, Figure~\ref{fig:chvatal_graph} shows two common ways to draw the Chv\'atal graph.

Users can also download all graphs matching their query (or a specific graph) in any of the supported formats (adjacency list, adjaceny matrix, graph6, and multicode) and then process these graphs further using their own software.

It is also possible to search for a specific graph by entering the HoG graph id of the graph, a graph6 string, or by drawing it using the graph editor (see Section~\ref{subsect:drawing}). This can be useful if you came across a specific graph in your research and you would like to check if it already appeared as an ``interesting'' graph in other contexts.

\subsubsection{Uploading new graphs}
\label{subsect:upload}

As already mentioned in Section~\ref{sect:introduction}, when the original House of Graphs website was launched in 2012, the searchable database was initialised with about 1\,500 graphs. But since then several new ``interesting'' graphs were added to the database and it currently contains nearly 22\,000 graphs.

To avoid spam or pollution of the database, only registered users can upload new graphs to the database (and for similar reasons only registered users can add comments or new graph drawings of existing graphs in the database). Any interested researcher can register a (free) account.

A registered user can upload a new graph either by drawing it using the graph editor (see Section~\ref{subsect:drawing}) or by uploading it from a file in one of the supported formats (e.g.\ graph6~\cite{graph6} or as adjacency matrix). Before the graph is added, the system first computes a canonical form of the graph using \textit{nauty}~\cite{nauty-website, MP14} to verify if the graph was not already present in the database. If it was already present, the user is notified about this and is redirected to the graph detail page of the existing graph. The user could then add a comment to the graph explaining that (s)he also came across this graph in the context of a different problem.

\begin{figure}[h!tb]
	\centering
	\includegraphics[width=0.8\textwidth]{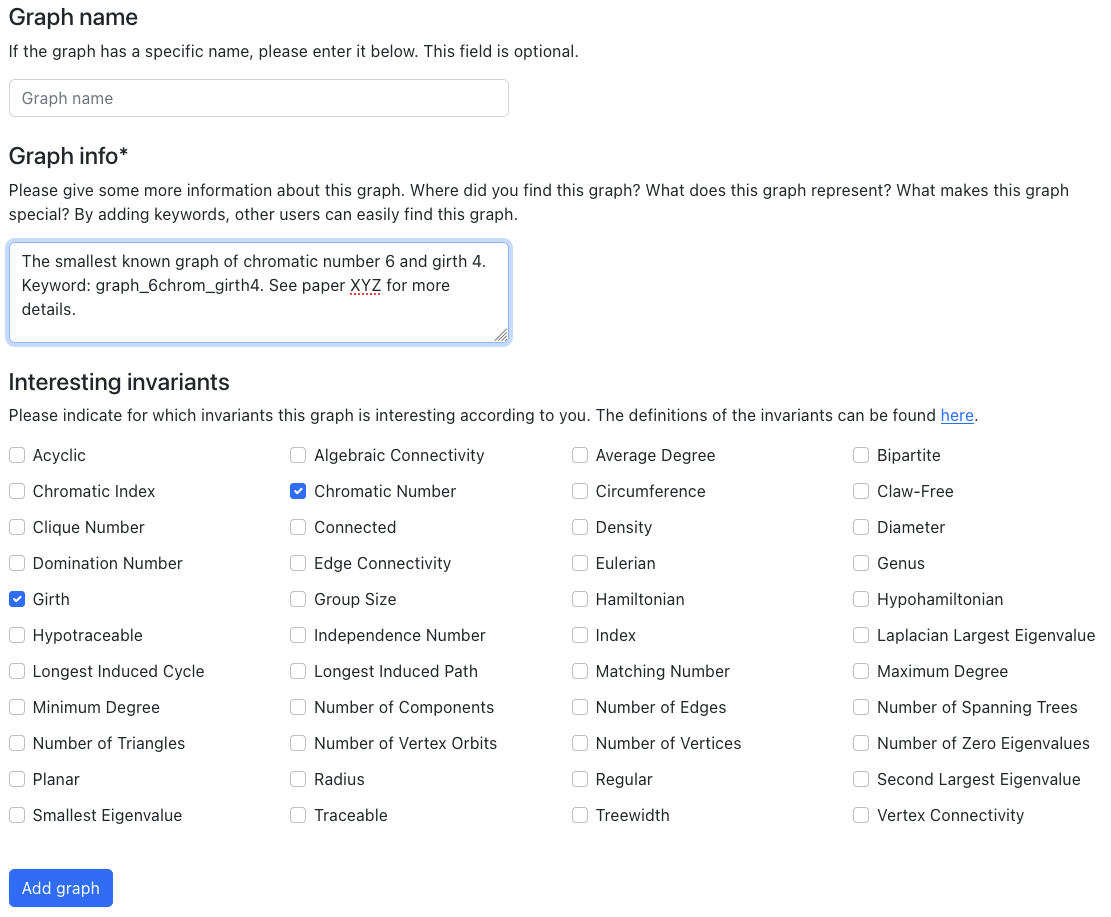}
	\caption{Screenshot of the page to add a new graph to the database.}
	\label{fig:add_graph}
\end{figure}

If the graph was not yet present in the database, the user is asked to enter some comments about this graph. For example: in what context this graph was obtained or some keywords which will allow others to find this graph more easily (so one can refer to this graph from a paper using this keyword). Optionally, the user can give the graph a name or mark the invariants for which (s)he believes the graph is especially interesting or relevant. (E.g.\ if the graph is the smallest graph of a given chromatic number and girth, it would make sense to mark ``chromatic number'' and ``girth'' as interesting invariants for this graph). Figure~\ref{fig:add_graph} shows a screenshot of how this is done on the website.

After the user has entered the relevant information about the graph, the graph is added to the database and a unique HoG graph id is generated so one can easily refer to this graph. In the background, jobs are submitted on the server to compute the invariant values of the new graph. As we support more than 40 invariants (see Table~\ref{table:invariants} for the complete list) and several of these invariants are known to be NP-hard to compute, we do not start these computations all at once and also do not let them run indefinitely.

\begin{figure}[h!tb]
	\centering
\includegraphics[width=0.4\textwidth]{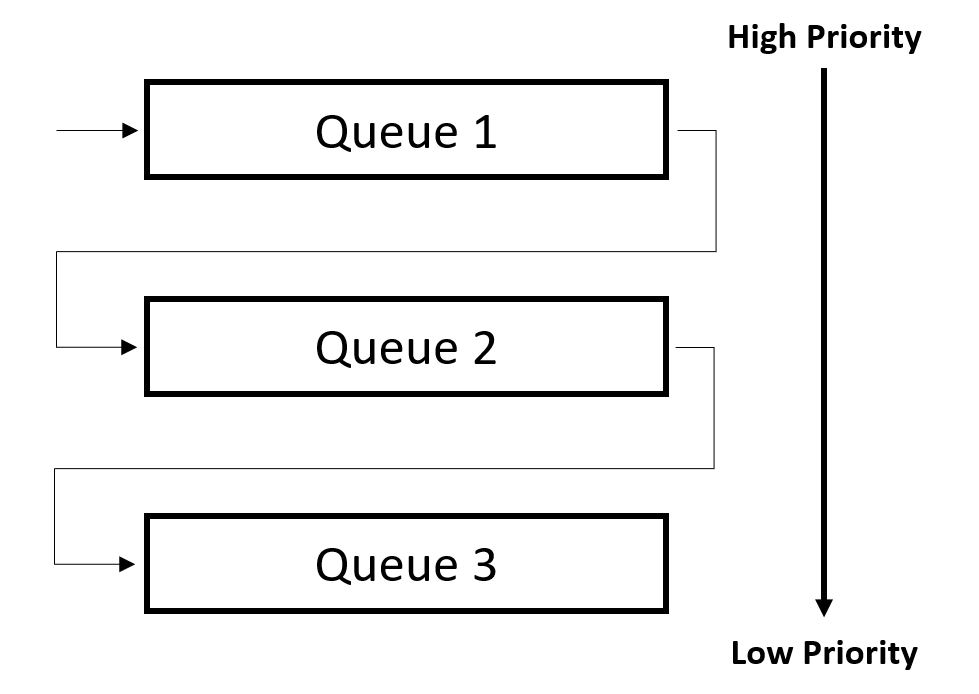}
	\caption{An example of a multilevel feedback queue.}
	\label{fig:multilevelqueue}
\end{figure}

Performing the calculations in a first-come-first-serve way would mean that fast processes are delayed by calculations that take a long time to complete. Our aim, however, is to calculate as many invariants as possible in the shortest amount of time, so the user can view a significant portion of the invariant values very shortly after uploading a graph. An additional problem arises when we look at the situation where multiple users upload a (different) graph almost simultaneously. In this case, we do not want to wait for all calculations on the first graph to finish before starting the calculations on the second graph. To overcome these problems, we decided to use a \textit{multilevel feedback queue}. Such a queue consists of multiple queues that can be used to handle calculations of different durations. Initially, all processes are assigned to the first queue, which we configured with a maximum execution time of one minute. Calculations that take longer than one minute to complete, are moved to the next queue, which has a longer execution time. This process can be repeated multiple times. By assigning a higher priority to queues with a smaller execution time, a system is created in which processes with small execution time are prioritised. An example of such a multilevel feedback queue can be seen in Figure~\ref{fig:multilevelqueue}.

We opted for the \textit{Slurm Workload Manager}~\cite{slurm} framework for the job scheduling. This is an open source framework which supports a multilevel feedback queue setup where priorities can be assigned to the queues. For more details about the queuing system we refer to~\cite{thesissven}.

When the graph is added (and while the invariant computations are still running or queued on the server), the user is immediately redirected to the graph detail page where (s)he can inspect the invariant values of this graph. The invariant values which were already computed are shown and invariants which still have to be computed or are currently being computed or could not be computed within the time limit of the final queue are marked as ``pending'',``computing'', or ``computation timeout'', respectively. On that page it is also possible to add alternative drawings of the graph (see Section~\ref{subsect:drawing}) and add/edit comments on the graph.

\subsection{Drawing graphs}
\label{subsect:drawing}

Another major feature of the House of Graphs which was significantly improved in the new website is its graph drawing and graph visualisation functionality. The old House of Graphs website already offered a very basic graph drawing tool, but it had some significant limitations. For example: it was not possible to zoom in/out (which can be very useful when handling large graphs), it was cumbersome to move vertices, there was no option to show the labels of the vertices, it was not possible to export the drawings, it was not possible to apply a spring embedder to an existing drawing, it was not possible to load a graph from a file (e.g.\ encoded as adjacency matrix), etc. Another major issue was that the old graph editor used the \textit{Google Web Toolkit} (\textit{GWT}) framework~\cite{gwt}, but this framework is no longer actively being developed as there were only three new updates in the past five years.

We therefore developed a new graph drawing tool which does not have these limitations. The new tool was written in \textit{D3}~\cite{d3}, which is a JavaScript library for producing dynamic, interactive data visualisations in web browsers which interacts well with \textit{React} (i.e.\ the framework we use for the front-end, cf.\ Section~\ref{sect:technologies}). Another advantage of \textit{D3} is its functionality to capture and handle user events. This makes it easy to detect click events, drag events,... and handle them accordingly (which was not possible or very cumbersome in the old tool).

Vertices can now be added by clicking on empty space on the canvas. By consecutively clicking on two vertices, an edge is created between these vertices. To accommodate for a fast and easy interaction with the graph, vertices can be dragged to new positions at any time. It is also possible to zoom in on the graph and pan around. This way, large or dense graphs can be edited more easily. On top of this, \textit{D3} makes it possible to apply forces to the vertices and edges to automatically create a good-looking structure. These automatic forces can be turned on or off by the user.

Another possible way to use the graph drawing tool is to load a graph from a file (e.g.\ encoded as adjacency matrix or in graph6 format). The tool then generates an automated spring embedding of the graph which can then still be manually fine-tuned by the user (and optionally the user can also ask to display the labels of the vertices). This is a useful way to generate an initial drawing of a graph and then manually analyse the structure of the graph by rearranging some of its vertices.

The resulting graph drawing can then be exported in various formats: both visualisation formats (e.g.\ svg, png, pdf, tikz,...) as well as graph encoding formats (as adjacency list, adjacency matrix, graph6 string or in multicode format). 

Users can also search if this graph was already present in the database or -- if the user is logged in -- add it to the database. In the latter case the drawing of the graph is also stored in the database and additional drawings of this graph can be added later (either by the same user or by other registered users).

\begin{figure}[h!tb]
	\centering
\begin{minipage}[c]{.52\textwidth}
\centering   
   \subfloat[]{\label{fig:graph_drawing}\includegraphics[width=1.0\textwidth]{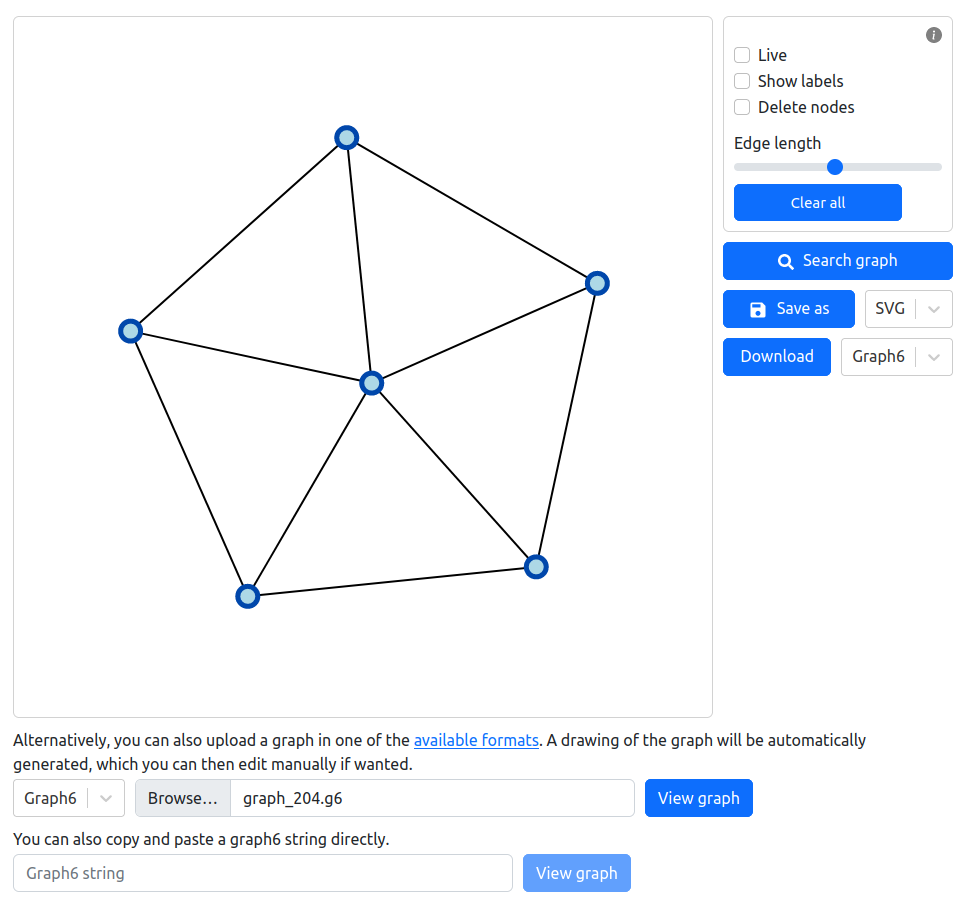}} 
\end{minipage}     	
\begin{minipage}[c]{.45\textwidth}
\centering   
   \subfloat[]{\label{fig:graph_drawing_old}\includegraphics[width=0.85\textwidth]{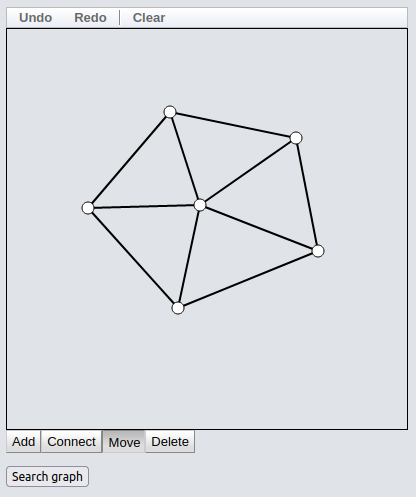}} 
\end{minipage}  	
	\caption{Screenshot of the new and old graph drawing tool, respectively.}
\end{figure}

Figure~\ref{fig:graph_drawing} shows a screenshot of the new graph drawing tool and for reference Figure~\ref{fig:graph_drawing_old} shows the old tool. 

\section{Concluding remarks}
\label{sect:conclusion}

As several of the frameworks used in the old House of Graphs website were approaching an end-of-life status while the website itself was still actively being used and cited by several researchers, it was in dire need of a major upgrade. With this complete rebuild of the House of Graphs website using modern frameworks, we believe that we have built a web application that is future-proof and much easier to maintain and extend. Moreover, we have added or improved several functionalities which improve the user experience such as the new graph drawing tool and the decoupling of the graph invariant computers which makes it much easier to add new graph invariants.

If readers of this article would have practical and well-tested programs to compute an interesting graph invariant which is not on the House of Graphs yet and which they would like to share, they are most welcome to contact us. Likewise if any readers would have a complete list of graphs which they would like to see added to the graph meta-directory. We also encourage interested readers to contribute to the website by uploading new graphs, adding comments and new graph drawings.

Next to adding more invariants and more graph classes, we still have several improvements or new functionalities in mind which we would like to add at some point in the future. This includes:

\begin{itemize}
\item Searching for graphs that satisfy some formula expressed in terms of other invariants, e.g.: $|Aut(G)| \geq |V(G)|$ (i.e.\ graphs where the order of their automorphism group is at least as large as their number of vertices).
\item Searching for graphs which (do not) contain a given graph as (induced) subgraph.
\item Linking graphs with their complement or line graph (if they are also present in the database).
\end{itemize}


\section*{Acknowledgements}
We would like to thank Gunnar Brinkmann, Brendan McKay and Jarne Renders who kindly provided programs to compute one or more graph invariants on the House of Graphs.
Next to that, we also thank the many users of the House of Graphs who uploaded new interesting graphs, added useful comments about graphs or contributed in any other way.



\bibliographystyle{plain}
\bibliography{references.bib}

\end{document}